\documentclass[12pt]{amsart}
\usepackage{amscd}
\usepackage{verbatim}
\usepackage{amssymb, amsmath, amsthm, amscd,ifthen}
\usepackage[dvips]{graphics}
\usepackage[cp866]{inputenc}
\usepackage{graphicx}
\usepackage{epsfig}
\usepackage{color}

\usepackage{amsmath,amssymb,amscd,}
\usepackage[all]{xy}

\usepackage[dvips]{graphics}
\usepackage[cp866]{inputenc}
\usepackage{graphicx}
\usepackage{epsfig}
\usepackage[hang]{subfigure}

\theoremstyle{plain}
\newtheorem{theorem}{Theorem}
\newtheorem{lemma}{Lemma}
\newtheorem{corollary}{Corollary}
\newtheorem{proposition}{Proposition}

\theoremstyle{definition}

\theoremstyle{plain}
\newtoks\thehProclaim
\newtheorem*{Proclaim}{\the\thehProclaim}

\theoremstyle{definition}
\newtoks{\thehRemark}
\newtheorem*{Remark}{\the\thehRemark}

\begin{document}

\title{A universality theorem for stressable graphs in the plane}

\author{ Gaiane Panina}

\address{G. Panina: St. Petersburg Department of Steklov Mathematical Institute, St. Petersburg State University, gaiane-panina@rambler.ru}

\keywords{Maxwell-Cremona correspondence, Grassmanian stratification, oriented matroid, equilibrium stress}

\begin{abstract}
  Universality theorems (in the sense of N. Mn\"{e}v) claim that the realization space of a combinatorial object (a point configuration, a hyperplane arrangement,  a convex polytope, etc.) can be arbitrarily complicated.
  In the paper, we prove a universality theorem for a graph in the plane with  a prescribed \textit{oriented matroid of stresses}, that is the collection of signs of all possible equilibrium stresses of the graph.

  This research is motivated by the Grassmanian stratification (Gelfand, Goresky,  MacPherson, Serganova) by thin Schubert cells, and by a recent series of papers on stratifications of configuration spaces of tensegrities (Doray, Karpenkov,  Schepers,  Servatius).
\end{abstract}

\maketitle

\section{Preliminaries and the main theorem}\label{section_preliminaries}
 Let  $\Gamma=(V,E)$ be a graph without loops and multiple edges, where $V= \{v_1,...,v_m\}$ is the set of vertices, and $E$ is the set of edges.
A \textit{realization} of  $\Gamma$ is a map $p:V\rightarrow \mathbb{R}^2$  such that $(ij)\in E$ implies $p(v_i)\neq p(v_j)$.
We abbreviate $p(v_i)$ as $p_i$.

That is, we have a planar drawing of $\Gamma$ with possible intersections of edges and possible coinciding vertices. However, each edge is mapped to a non-degenerate line segment.

A \textit{stress} $\mathfrak{s}$ on a realization $(\Gamma, p)$ is an assignment of real scalars  $\mathfrak{s}(i,j)$  to the edges.
One imagines that each edge is represented by a (either compressed or extended) spring. Each spring produces some forces at its endpoints.

A stress $\mathfrak{s}$ is called a \textit{self-stress}, or an \textit{equilibrium stress}, if at every vertex  $p_i$,  the sum of the forces produced by the springs vanishes:
$$\sum_{(ij)\in E}\mathfrak{s}(i,j)\mathbf{u}_{ij}=0.$$
Here $\mathbf{u}_{ij}=\frac{p_i-p_j}{|p_i-p_j|}$ is the unit vector pointing from $p_j$ to $p_i$.

A self-stress is \textit{non-trivial} if it is not identically zero.

The set of all self-stresses $\mathfrak{S}(\Gamma,p)$ is a linear space which naturally embeds in $\mathbb{R}^e$, where $e=|E|$; the space $\mathfrak{S}$ depends on $p$.

A realization $(\Gamma,p)$ is  \textit{ stressable} if  $dim~\mathfrak{S}(\Gamma,p)>0.$

Given  $(\Gamma,p)$, define an oriented matroid $\mathcal{M}(\Gamma,p):=SIGN(\mathfrak{S}(\Gamma,p))$.

In simple words,  to obtain the matroid, enumerate somehow the edges of the graph, and for each non-trivial stress, list the signs of $\mathfrak{s}(i,j)$.
We obtain a collection of strings (elements of $(+,-,0)^{\sharp (E)}$), which is an oriented matroid\footnote{This is some realizable oriented matroid indeed, since it represents the set of vectors of some easy-to-build vector configuration related to the \textit{rigidity matrix}. For rigidity matrices see \cite{W}.}. For the purpose of the present paper, it is sufficient to imagine an oriented matroid as a collection of strings. For a general theory of oriented matroids see \cite{OM}.

\medskip

\textbf{Example:} let $(\Gamma, p)$ be a planar realization of the graph $K_4$  such that $p_4$ lies inside the triangle $p_1p_2p_3$.
Assume that the edges are enumerated in a way such that first come the edges of the triangle.
Then $\mathcal{M}(\Gamma,p)=\{(+++---),(---+++)\}$.

\medskip

 The \textit{realization space } of a graph $\Gamma$
is the space of all realizations of $\Gamma$  factorized   by the action of the affine group:

 $$\mathcal{R}(\Gamma)=\{p:  p \hbox{ is a realization of } \Gamma\}/Aff(\mathbb{R}^2).$$

Given a graph $\Gamma$ and an oriented matroid $\mathcal{M}$,  the \textit{realization space } of $(\Gamma,\mathcal{M})$
is the space of all realizations of $\Gamma$ that yield the oriented matroid $\mathcal{M}$:

 $$\mathcal{R}(\Gamma, \mathcal{M})=\{p\in \mathcal{R}(\Gamma): \mathcal{M}(\Gamma,p)= \mathcal{M}\}.$$

For a fixed graph $\Gamma$, the realization spaces $\mathcal{R}(\Gamma, \mathcal{M})$ stratify  $\mathcal{R}(\Gamma)$.
Each of $\mathcal{R}(\Gamma, \mathcal{M})$ becomes a \textit{stratum}.

\medskip

In general, semialgebraic sets are subsets of some Euclidean space  $\mathbb{R}^N$ defined by
polynomial equations and inequalities.
A  semialgebraic set is called a \textit{open basic primary semialgebraic set (OBP semialgebraic set)} if there are no defining equations,
all the defining inequalities are strict, and the coefficients of all the defining polynomials are rational.

We borrow the notion of stable equivalency from traditional papers on universality, e.g. from \cite{RG}: \textit{stable equivalence} is an equivalence relation on OBP semialgebraic sets generated by
rational equivalence and stable projections.

\bigskip
The main result of the paper is:
\begin{theorem}\label{MainTheorem}For each open basic primary  semialgebraic set $\mathfrak{A}$, there exists a graph $\Gamma$ and an oriented matroid $\mathcal{M}$
such that the realization space $\mathcal{R}(\Gamma, \mathcal{M})$ is stably equivalent to $\mathfrak{A}$.
\end{theorem}

Our first motivation comes from the complex Grassmanian stratifications \cite{GGMS}, where strata are labeled by   realizable matroids, and  each stratum equals the realization space of a matroid.  The stratification has a version over the field $\mathbb{R}$, where  strata are labeled by realizable oriented matroids.

 The other motivation is  a series of  papers \cite{DorKarp}, \cite{Karp}, \cite{KarpSS} on stratifications of configuration spaces of tensegrities. Although the setup of the present paper might look different from the setup of  \cite{DorKarp}, \cite{Karp}, \cite{KarpSS}, there is very much in common, see Section 3. In  particular, we have the following  analogue of Theorem \ref{MainTheorem}:

\begin{theorem}\label{MainTheorem2}
For each connected open basic primary  semialgebraic set $\mathfrak{A}$, there exists a graph $\Gamma$ and a stratum (in the sense of  \cite{DorKarp}) $R$ in the realization space of $\Gamma$ such that $R$ is stably equivalent to $\mathfrak{A}$.
\end{theorem}

\medskip

\textbf{Acknowledgement.} The author is grateful to Joseph Gordon and Yana Teplitskaya for multiple discussions.  This research is supported by the Russian Science Foundation under grant 16-11-10039.

\section{Proof of Theorem \ref{MainTheorem}}\label{}
Combinatorial equivalence of planar point configurations and line configurations see \cite{RG2} is a classical subject  and a starting point of our research.
Given a combinatorial type, the equivalence class of point configurations (or line configurations) having this type is called the \textit{realization spaces}.
We shall use the same letter $\mathcal{R}$ for realization spaces of line configurations. In particular, for a line configuration $L$ we denote
by $\mathcal{R}(L)$ the realization space of all line configurations that are combinatorially equivalent to $L$.

We shall use the following version (chronologically, one of the first ones) of the celebrated Universality Theorem \cite{Mnev}: \textit{generic  planar point configurations are universal}. More precisely,
for each OBP semialgebraic set $\mathfrak{A}$, there exists a planar point configuration with points in generic position\footnote{No three points are collinear.} such that the realization space of the configuration is stably equivalent to $\mathfrak{A}$. An immediate consequence of the theorem is: generic planar line arrangements are universal.

\begin{figure}
\centering
\includegraphics[width=12 cm]{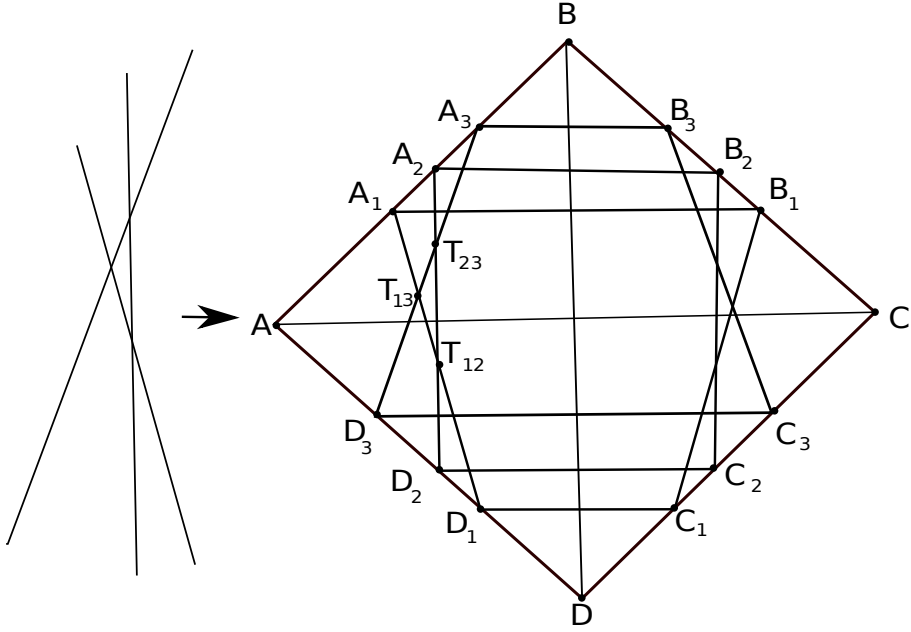}
\caption{The graph $\Gamma$ with its realization  $p$.  }\label{Figure_MainGraph}
\end{figure}

Assume that a  OBP semialgebraic set $\mathfrak{A}$  is fixed.
For the set $\mathfrak{A}$, we shall construct a graph $\Gamma$ together with its realization $p$, depicted in Fig. \ref{Figure_MainGraph}. Thus, we get the associated oriented  matroid $\mathcal{M}=\mathcal{M}(\Gamma,p)$.  Our final aim is to show that $\mathcal{R}(\Gamma, \mathcal{M})$ is stably equivalent to $\mathfrak{A}$.

\medskip

\medskip

Here is the construction.

(1) Take   a generic line configuration $L=\{l_i\}_{i=1}^n$ whose realization space is stably equivalent to $\mathfrak{A}$. Since the configuration is generic, there are no triple intersections.

(2) Take a rhombus $ABCD$ such that all mutual intersections $T_{ij}=T_{ji}=l_i \cap l_j$ lie strictly inside the rhombus, and each line $l_i \in L$ intersects the interiors both of the segments $AB$ and $AD$. Denote the intersection  points points by $A_i$ and $D_i$ respectively. We may assume that the points $A,A_1,A_2,...,A_n,B$ appear on the segment $AB$ in this very order. Therefore, the order of the points $D_i$ (from left to right) is reverse.

(3) Add to our construction the diagonals of the rhombus $AC$ and $BD$.

(4) Add to our construction the points $B_i\in BC$, $C_i\in CD$, and the segments $A_iB_i$, $B_iC_i$, $C_iD_i$, and $D_iA_i$ such that
$A_iD_i$ is symmetric to $B_iC_i$ with respect to the diagonal $BD$ for all $i$.

(5) Finally, add the intersection points $T_{ij}=T_{ji}=A_iD_i\cap A_jD_j$.

(6) Now let us specify edges of the graph. The points $A_i$ split the segment $AB$ into edges. The points $B_i$ split $BC$ into edges, etc. Besides, the points $T_{ij}$
split $A_iD_i$ into edges. The segments $A_iB_i$, $B_iC_i$, $C_iD_i$,  $AC$, and $BD$ are edges as well. All the edges are depicted in Fig. \ref{Figure_MainGraph}.

We obtain a realization of a graph, whose vertices are \newline $\{A,B,C,D,A_i,B_i,C_i,D_i,T_{ij}\}_{i,j}$.

\begin{figure}
\centering
\includegraphics[width=6 cm]{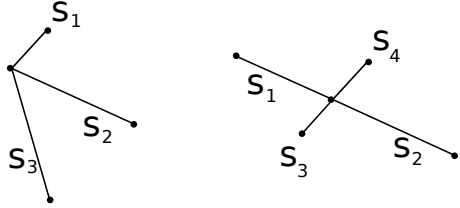}
\caption{Local stresses }\label{Figure_Local1}
\end{figure}

\begin{lemma}\label{LemmaLocal}\begin{enumerate}
                                 \item Assume that we have a self-stressed realization of an arbitrary graph, and a vertex  of the graph looks as is depicted in Fig. \ref{Figure_Local1}.
Then, in notation of the figure,  the  stresses satisfy:\begin{enumerate}
                                                                        \item (left) $Sign~s_1=Sign~s_3=-Sign~s_2$,
                                                                        \item (right) $s_1=s_2, \ \ s_3=s_4$.
                                 \item If a self-stress $\mathfrak{s}$ of the above constructed $(\Gamma, p)$ (see Fig. \ref{Figure_MainGraph}) vanishes at one of the segments of the quadrilateral $A_iB_iC_iD_i$, then it vanishes on each of the segments of the quadrilateral.
                                \item    If a self-stress $\mathfrak{s}$ of $(\Gamma, p)$  from Fig. \ref{Figure_MainGraph} vanishes at all  the segments lying on  $AB$, then it vanishes everywhere.
                                     \qed
                               \end{enumerate}
            \end{enumerate}
\end{lemma}

Let us look at some particular elements of $\mathfrak{S}(\Gamma,p)$ (in matroid terminology, they all are \textit{circuits} of the oriented  matroid $\mathcal{M}$). At most of the edges,  these stresses vanish, so
we depict  them as  subgraphs of $(\Gamma,p)$. That is, we leave stressed edges only, and indicate the signs of the stress.

\begin{figure}
\centering
\includegraphics[width=14 cm]{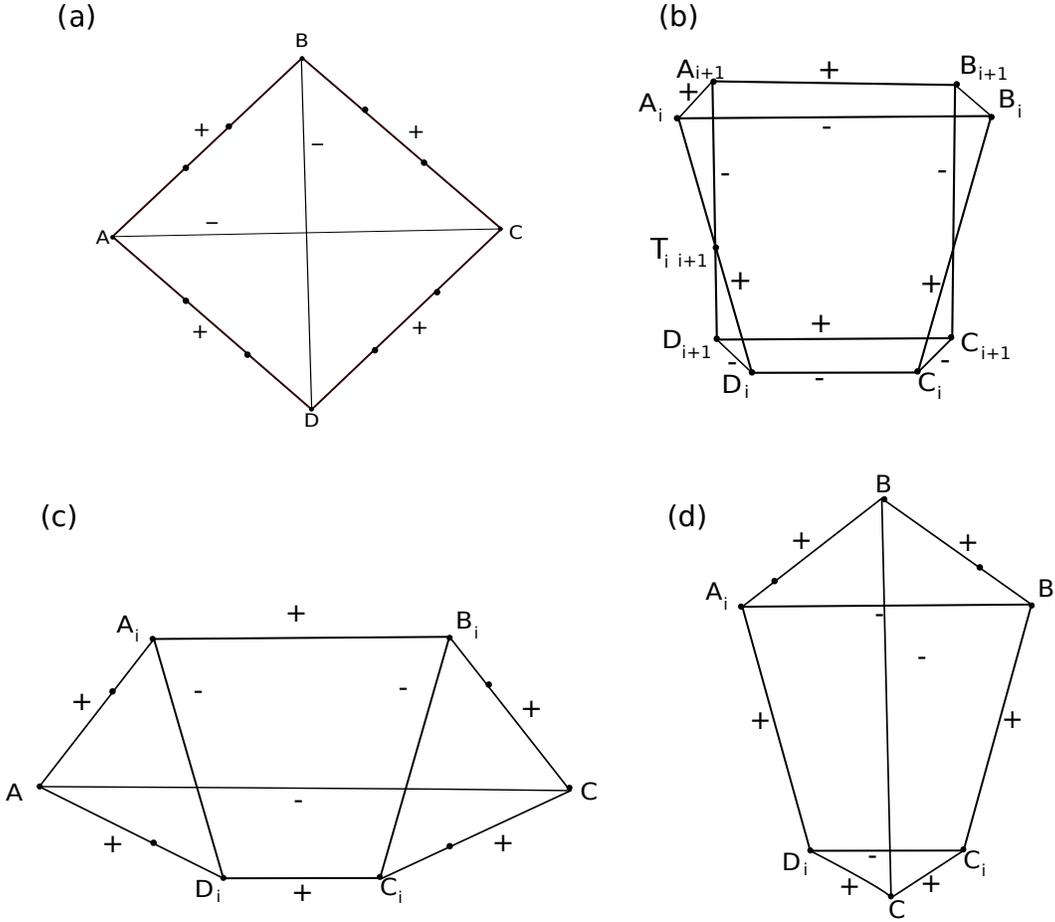}
\caption{Some particular stresses.  }\label{Figure_Circ}
\end{figure}

\begin{lemma}\label{LemmaCirc}\begin{enumerate}
                                \item The subgraphs depicted in Figure \ref{Figure_Circ}  are stressable. The signs of the associated stresses are indicated. (Clearly, simultaneous inversion of signs also represents some self-stress.)

                                \item The stresses (a) and (d) from Figure  \ref{Figure_Circ} (for all $i=1,...,n$) form a basis of the linear space $\mathfrak{S}(\Gamma,p)$.
                                 \item The stresses (a) and (c) (for all $i=1,...,n$) also form a basis of $\mathfrak{S}(\Gamma,p)$.
                                 \item For any stress of $(\Gamma,p)$,   the ratio of stresses on edges $A_iA_{i+1}$ and $B_iB_{i+1}$ does not depend on $i$ (provided that the stresses are non-zero).
                                    The ratio of stresses on edges $C_iC_{i+1}$ and $D_iD_{i+1}$ does not depend on $i$ either.
                              \end{enumerate}

\end{lemma}
Proof. (1) (a) is known to be stressable.\footnote{(a) can be viewed as a projection of a tetrahedron, therefore (a) is \textit{liftable}. By Maxwell-Cremona correspondence, it is stressable.}   (c) and (d) are stressable since these are \textit{Desargues configurations}. This means that the three lines $A_iD_i$, $BC$ and $B_iC_i$ meet at a point;  the three lines $A_iB_i$, $C_iD_i$ and $AC$ are parallel, that is, meet at a point at infinity.  (b) is the difference of two  different stresses of type (c), therefore stressable. The signs in all the cases follow from Lemma \ref{LemmaLocal}.

Prove (2). Assume we have a stress $\mathfrak{s}\in \mathfrak{S}(\Gamma,p)$.  Adding an appropriate stress of type (a), we kill the value of the stress on the edge $AA_1$, and therefore, on all the edges emanating from $A$. Next, adding  appropriate stresses of type (d)  kills the stresses on all the edges of $AB$. By Lemma \ref{LemmaLocal}, the result is identical zero.

(3) follows from (2).  (4) is true for (a) and (d), therefore, it is true for all self-stresses.\qed

\medskip
Now  we analyze the realization space of  matroid $(\Gamma, \mathcal{M})$.

\begin{proposition}\label{prop_1}
Assume that  $\mathcal{M}(\Gamma, p')=\mathcal{M}(\Gamma, p)$, that is, $(\Gamma, p')\in \mathcal{R}(\Gamma, p).$
Denote the vertices of the realization by the same letters with primes (that is, by $A'_i, B_i',$ etc).
Then \begin{enumerate}
       \item All collinearities of vertices that are present in $p$ survive for $p'$. Besides, the order of collinear vertices maintains.
       \item The points $A'B'C'D'$ lie in the convex position.
       \item $p'$ yields an arrangement of lines $L'=\{l_i'\}$ with the same combinatorics as the initial arrangement $L$.
     \end{enumerate}
\end{proposition}
Proof. The stressed graph (a) from Fig. \ref{Figure_Circ} should be  stressed with the same signs (and with no other signs)  for $p'$ as well. Therefore $A'B'C'D'$ is a convex quadrilateral, and  all the edges belonging to $A'B'$, are  collinear (all the edges belonging to $C'B'$, etc. are collinear as well).
The graphs of type (b) from Fig. \ref{Figure_Circ}  remain stressed for $p'$,  so the segments of $l_i$ stay collinear. Therefore we have an arrangement of lines $L'=\{l_i'\}$

The graph $\Gamma$ records the combinatorics of $L$, therefore the $L$ and $L'$ are combinatorially equivalent, that is, $L'\in \mathcal{R}(L)$.\qed

\begin{corollary}\label{corol}There exists a natural mapping between the realization space of $(\Gamma, \mathcal{M})$ and the realization space of the arrangements of lines $L$:
$$\pi: \mathcal{R}(\Gamma, \mathcal{M})\rightarrow \mathcal{R}(L).$$
The mapping $\pi$  extracts the arrangement $L'$ from $(\Gamma, p')$ and forgets the rest.\qed
\end{corollary}

\begin{proposition}\label{Prop2}
  Let $A'B'C'D'$ be a convex quadrilateral.
  Assume that the points  $\{A'_i,B'_i,C'_i,D'_i\}_{i=1}^n$ are such that
  \begin{enumerate}

 \item the points  $A'_i$ lie
  on the segment $A'B'$ and come in the same order as $A_i$. The points  $B'_i$ lie
   on the segment $B'C'$ and come in the same order as $B_i$;
      and  the same condition for $C'_i$ and $D'_i$.
     \item The affine hulls of $A_i'D'_i$ form an arrangement of lines combinatorially equivalent to $L$.
     \item All the associated subgraphs of type (c) from Fig. \ref{Figure_Circ} are Desargues ones.
     \end{enumerate}

     Then\begin{enumerate}
           \item All the associated subgraphs of type (d) from Fig. \ref{Figure_Circ} are Desargues ones, and therefore, stressable.
           \item The realization of the graph $\Gamma$ with these vertices has the same oriented matroid as $\mathcal{M}=\mathcal{M}(\Gamma, p)$.
         \end{enumerate}
\end{proposition}
Proof.
(1) follows from Desargues'  theorem. The conditions imply that we have some realization $p'$ of $\Gamma$, and that all circuits depicted in
Fig. \ref{Figure_Circ}
are circuits relative $p'$.

\begin{figure}
\centering
\includegraphics[width=4 cm]{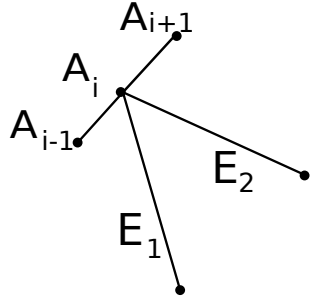}
\caption{Illustration for the proof of Proposition \ref{Prop2}  }\label{Figure_Local2}
\end{figure}

Before we proceed with the claim (2), let us observe the following:
\begin{lemma}\label{small lemma} Assume that $\mathfrak{s}$ is a self-stress of $(\Gamma,p')$, whose values on the edges  $A'_{i-1}A'_i$ and $A'_{i+1}A'_i$
we denote by $s_1$ and $s_2$. Then the signs of the stresses of the other two edges $E_1$ and $E_2$ (each of them equals $A'_iT'_{ij}$ for some $j$), emanating from $A_i'$ are:
$$SIGN(\mathfrak{s}(E_1))=SIGN(s_2-s_1)=- SIGN(\mathfrak{s}(E_2)),$$
assuming that $E_1$ lies to the left of $E_2$, see Fig. \ref{Figure_Local2}.
Similar statements are valid for edges emanating from $B'_i, C'_i$, and $D'_i$.
\qed
\end{lemma}

Now let us prove the statement (2) of Proposition \ref{Prop2}. First observe that Lemma \ref{LemmaCirc} stays valid for $(\Gamma,p')$.
Let $\mathfrak{d}_i$ (respectively, $\mathfrak{d}'_i$) be a stress of $(\Gamma,p)$ (respectively, $(\Gamma,p')$) depicted in Fig. \ref{Figure_Circ}, (d), such that its value on $AA_1$ (respectively, $A'A_1'$) equals $1$.

Let $\mathfrak{a}$ be a stress of $(\Gamma,p)$ depicted in Fig. \ref{Figure_Circ}, (a), such that its value on $AA_1$ equals $1$,
let $\mathfrak{a}'$ be defined analogously for $(\Gamma,p')$.

 Assume that $\mathfrak{s}$ is a stress of $(\Gamma,p)$. By Proposition \ref{prop_1}, $\mathfrak{s}=\lambda \mathfrak{a}+\sum \lambda_i\mathfrak{d}_i$ for some real coefficients.
 Consider the stress $\mathfrak{s}'=\lambda \mathfrak{a}'+\sum \lambda_i\mathfrak{d}'_i$ of $(\Gamma,p')$. By Lemma \ref{LemmaCirc} (4) and Lemma \ref{small lemma}, we have
  $SIGN(\mathfrak{s})=SIGN(\mathfrak{s}')$. Conversely, each stress $\mathfrak{s}'$ of $(\Gamma,p')$, has a similar counterpart for $(\Gamma,p)$.\qed

\begin{proposition}\label{lastProp} $\pi$ is a stable projection.
\end{proposition}

Proof. Assume that a line arrangement $L'$ belongs to the realization space $\mathcal{R}(L)$. Let us look at the preimage $\pi^{-1}(L')$. Specification of $A'B'C'D'$ is stable since the positions of the lines $A'B', \ A'D'$, etc. are defined by a number of strict inequalities depending on the lines $L$.

Now we may choose arbitrary distinct points $A_1',...,A'_n$ on the segment $A'B'$ that come in the same order as $A_1, A_2,...,A_n$.
The same happens with $B_1',...,B'_n$: here we care about their order only.
Now let us choose $C'_1$ as an arbitrary point on the segment $B'C'$. Once $C'$ is specified, the position of $D_1'$ is uniquely determined, since Desargues condition implies that the lines $A'C'$, $ A'_1B'_1$, and $C'_1D'_1$  meet at a point.

The point $C_2'$ should be chosen to the left of $C_1'$ but in such a way that $D_2'$ lies to the right of $D'_1$. This is always possible but dictates some extra condition,  still in the framework of stable equivalence.
The rest of the points $C'_i$ and $D'_i$ are treated analogously. \qed

\medskip

 Corollary \ref{corol} and Proposition \ref{lastProp} imply that  $\mathcal{R}(\Gamma, \mathcal{M})$ is stably equivalent to $\mathfrak{A}$. Theorem 1 is proven.

\section{Relations between different settings. Proof of Theorem \ref{MainTheorem2}.}

In the section we show the equivalence of  the settings of the present paper and that of  \cite{DorKarp}.

Let us start with the definition of equilibrium stress. The paper \cite{DorKarp} puts no restrictions on a realization of a graph $p$, that is, the endpoints of an edge might be mapped to one and the same point. Besides, \cite{DorKarp} presents
a more usual setting of the equilibrium stress (as  in \cite{Conn}):  the equilibrium condition reads as $$\sum_{(ij)\in E}{s}(i,j)(p_i-p_j)=0.$$

Let us denote by $S(\Gamma, p)$ the linear  space of stresses and by $M(\Gamma,p)=SIGN(S(\Gamma,p))$ the associated matroid.

Clearly, if no edge is degenerate (that is,  $p_i - p_j\neq 0$), a stress $s$ in this setting gives a stress in the setting of the present paper $\mathfrak{s}(i,j)=s(i,j)|p_i-p_j|$
and vice versa. Therefore, ${M}(\Gamma,p)=\mathcal{M}(\Gamma,p)$. The only subtlety may arise if a realization $p$  produces degenerate edges.

\begin{lemma}\label{lemmaDiffSett}
 The matroid ${M}(\Gamma,p)$ "knows" all the degenerate edges. In particular,  if there exists $p\in R(\Gamma, M)$ with no degenerate edges,
 then each $p'\in R(\Gamma, M)$ has no degenerate edges.

\end{lemma}
Proof. Degenerate edges are detected by almost everywhere zero stresses: an edge number $i$ is degenerate for $(\Gamma, p)$ iff   $({0,...,0,+,0,0,...,0})\in M(\Gamma,p)$.\qed

\subsection*{Strong equivalence vs weak equivalence}  Assume that a realization $p$ of a graph $\Gamma$ has no degenerate edges.

Repeating \cite{KarpSS}, let us say that two realizations of one and the same graph $(\Gamma,p)$ and $(\Gamma,p')$ are \textit{strongly equivalent},
if there exists a sign preserving homeomorphism between the stress spaces $\mathfrak{S}(\Gamma,p)$ and $\mathfrak{S}(\Gamma,p')$.

Two realizations of one and the same graph $(\Gamma,p)$ and $(\Gamma,p')$ are \textit{weakly equivalent},
if the associated  matroids coincide: $\mathcal{M}(\Gamma,p)=\mathcal{M}(\Gamma,p')$.

Classes of weak equivalence  are realization spaces, defined in the Introduction.
Classes of strong equivalence are strata considered in \cite{KarpSS}.\footnote{To be more precise, in \cite{KarpSS}  the strata are the connected components of the classes of strong equivalence.}

\begin{proposition}\label{PropFine}
Strong equivalence  equals weak equivalence.
\end{proposition}

Proof. Clearly, strong equivalence implies weak equivalence. Let us prove the converse.
The linear space $\mathfrak{S}(\Gamma,p)$  is tiled by convex cones, each  cone corresponds to some string of signs from
$\mathcal{M}(\Gamma,p)$.  Let us intersect this tiling with the unit sphere centered at the origin. This gives a  tiling of the sphere  where each tile is a  spherically convex polytope. The matroid $\mathcal{M}(\Gamma,p)$ "knows" the incidence relation of the tiles: a tile labeled by
$(\varepsilon_1,...,\varepsilon_e)$, $\varepsilon_i \in \{+,-,0\}$ belongs to the closure of the tile $(\varepsilon'_1,...,\varepsilon'_e)$
iff either  $\varepsilon_i=0$, or $\varepsilon_i=\varepsilon_i'$.

Besides, the matroid $\mathcal{M}(\Gamma,p)$ "knows" the dimension of each tile.
Indeed, the matroid records the face poset of each tile. Since each tile is  some pointed cone, its dimension is determined by the length of a longest chain in the poset.

  Now it becomes possible to inductively build a sign-preserving homeomorphism between two spaces $\mathfrak{S}(\Gamma,p)$ and $\mathfrak{S}(\Gamma,p')$ with equal matroids.
One should  start with zero-dimensional spherical tiles, then extend the homeomorphism to one-dimensional tiles, etc.
\qed

\medskip

Now let us prove Theorem \ref{MainTheorem2}. Given  $\mathfrak{A}$, take the pair $(\Gamma,p)$ as in the proof in Theorem \ref{MainTheorem}.
By Lemma \ref{lemmaDiffSett}, $R(\Gamma,p)=\mathcal{R}(\Gamma,p)$. By Proposition \ref{PropFine}, $\mathcal{R}(\Gamma,p)$ is a strong equivalence class, that is, a stratum in the sense of  \cite{DorKarp}, \cite{Karp}, \cite{KarpSS}. Finally, by Theorem \ref{MainTheorem} the stratum is stably equivalent to $\mathfrak{A}$.

\section{Appendix}

\subsection*{One more example} A simpler (but in a sense, more "degenerate") example of   $(\Gamma', p)$ with the same realization space as in the previous section can be obtained if one takes $(\Gamma,p)$ from
Fig. \ref{Figure_MainGraph}, removes all the edges lying on $AB$, $BC$, $CD$, $DA$, $AC$, $BD$ $A_iB_i$, $B_iC_i$, $C_iD_i$, and adds the new edges $A_iD_i$. That is, all the edges of the new graph lie on the lines $l_i$.
\subsection*{"No parallel edges" condition} 
\begin{figure}
\centering
\includegraphics[width=6cm]{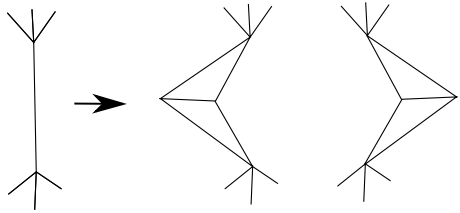}
\caption{ Adding a stressed copy of $K_4$.  }\label{Figure_non-par}
\end{figure}

The graph from Figure \ref{Figure_MainGraph}  has  parallel edges emanating from one and the same vertex (in fact, almost all the vertices have parallel emanating edges).
If  parallel edges emanating of one and the same vertex are forbidden we still have a universality-type theorem:

\begin{theorem}\label{MainThm}For each OBP semialgebraic set $\mathfrak{A}$, there exists a graph $\Gamma$, an oriented matroid $\mathcal{M}$, and a number $N$
such that \begin{enumerate}
            \item $(\Gamma, \mathcal{M})$ has a realization with no parallel edges at all, and
            \item the realization space $\mathcal{R}(\Gamma, \mathcal{M})$ is stably equivalent to $2^N$ disjoint copies of $\mathfrak{A}$.
          \end{enumerate}

\end{theorem}
Proof. The idea is depicted in Fig. \ref{Figure_non-par}: take the graph from Figure \ref{Figure_MainGraph} and for each edge $(ij)$, add two new vertices and replace $(ij)$  by five new edges.  One imagines a stressed realization of $K_4$ added to a stressed realization of $\Gamma$ in such a way that the stresses
on $(ij)$ cancel. Denote the realization of the new graph by $(\hat{\Gamma},\hat{p})$, and set $\hat{M}=M(\hat{\Gamma},\hat{p})$.
There exists a natural mapping $$R(\hat{\Gamma},\hat{M})\rightarrow R({\Gamma},{M}).$$

The preimage of each point has $2^{e(\Gamma)}$ connected components since
 each stressed
 copy of $K_4$ can be attached both on the righthand side and on the lefthand side of $(ij)$, but never degenerates. \qed

\subsection{Intersection of closures of two strata is not necessarily the  closure of a stratum}  This phenomenon was observed in \cite{GGMS} for Grassmanian stratifications. Let us adjust an example from  \cite{GGMS} to show the same for stressed graphs.

\begin{figure}
\centering
\includegraphics[width=7cm]{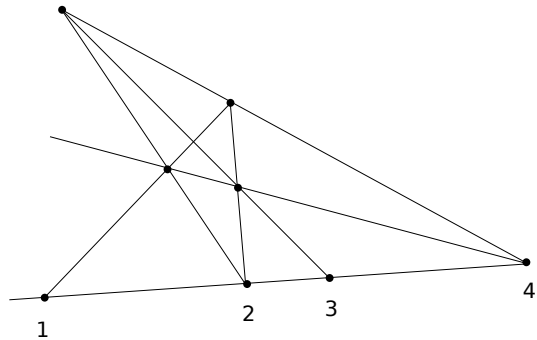}
\caption{A graph which forces harmonic relation.  }\label{Figure_stratum}
\end{figure}

Take the point configuration from  Fig. \ref{Figure_stratum} and associate to it a graph $(\Gamma,p)$ by the following rule: for each three collinear points $i,j,k$ add the edges $(ij)$, $(jk)$, and $(ik)$. So each three collinear points yield a stressable subgraph $K_3$.  We conclude that all the collinearities of vertices persist for all the elements of the realization space   $\mathcal{R}(\Gamma,\mathcal{M}(\Gamma,p))$.  However, all these collinearities imply that the four points $1,2,3,4$ are harmonic, that is, their cross ratio equals $-1$.

Let us take a realization $p'$ of $\Gamma$ with all the vertices lying on a line. The corresponding matroid depends on the order of the vertices only and "does not see" the cross ratio. Therefore the intersection of the closures of the strata $\mathcal{R}(\Gamma,\mathcal{M}(\Gamma,p))$ and $\mathcal{R}(\Gamma,\mathcal{M}(\Gamma,p'))$ is not a closure of a  stratum.

\end{document}